\documentclass[letterpaper, 10 pt, conference]{ieeeconf}  

\IEEEoverridecommandlockouts                              

\def\BibTeX{{\rm B\kern-.05em{\sc i\kern-.025em b}\kern-.08em
    T\kern-.1667em\lower.7ex\hbox{E}\kern-.125emX}}

\usepackage{amsthm}
\usepackage{cite}
\usepackage{amsmath,amssymb,amsfonts}

\usepackage{algorithmic}
\usepackage{graphicx}
\usepackage{textcomp}
\usepackage{xcolor}
\usepackage{verbatim}
\usepackage{caption}
\usepackage{subcaption}
\usepackage{tabularray}
\usepackage{multirow}
\allowdisplaybreaks
\usepackage{etoolbox}
\usepackage{epstopdf}

\newtheoremstyle{mystl}
  {0} 
  {0} 
  {\itshape}
  {}
  {\bfseries}  
  {.} 
  {5pt plus 1pt minus 1pt} 
 {\thmname{#1} \thmnumber{#2}\ifblank{#3}{}{ (\thmnote{#3})}} 
\theoremstyle{mystl}

\theoremstyle{mystl}

\newtheorem{thme}{Theorem}

\newtheorem{ass}{Assumption}

\newtheorem{rmk}{Remark}
\usepackage{mathtools}

\begin{document}
\title{\LARGE \bf
Setpoint Tracking and Disturbance Attenuation for Gas Pipeline Flow Subject to Uncertainties using Backstepping
}

\author{Bhathiya Rathnayake$^{1}$, Anatoly Zlotnik$^{2}$, Svetlana Tokareva$^2$, and Mamadou Diagne$^{3}$
\thanks{$^{1}$B. Rathnayake is with the Department of Electrical and Computer Engineering, University of California San Diego, USA,  
        {\tt\small brm222@ucsd.edu}}%
\thanks{$^{2}$A. Zlotnik and S. Tokareva are with the Applied Mathematics and Plasma Physics Group, Los Alamos National Laboratory, USA}%
        \thanks{$^{3}$M. Diagne is with the Department of Mechanical and Aerospace Engineering, University of California San Diego, USA }%
}

\maketitle



\begin{abstract}
In this paper, we consider the problem of regulating the outlet pressure of gas flowing through a pipeline subject to uncertain and variable outlet flow. Gas flow through a pipe is modeled using the coupled isothermal Euler equations, with the Darcy–Weisbach friction model used to account for the loss of gas flow momentum. The outlet flow variation is generated by a periodic linear dynamic system, which we use as a model of load fluctuations caused by varying customer demands. We first linearize the nonlinear equations around the equilibrium point and obtain a $2 \times 2$ coupled hyperbolic partial differential equation (PDE) system expressed in canonical form. Using an observer-based PDE backstepping controller, we demonstrate that the inlet pressure can be manipulated to regulate the outlet pressure to a setpoint, thus compensating for fluctuations in the outlet flow. Furthermore, we extend the observer-based controller to the case when the outlet flow variation is uncertain within a bounded set. In this case, the controller is also capable of regulating the outlet pressure to a neighborhood of the setpoint by manipulating the inlet pressure, even in the presence of uncertain fluctuations in the outlet flow. We provide numerical simulations to demonstrate the performance of the controller.
\end{abstract}

\section{Introduction}
In recent years, significant attention has been directed to study the system dynamics in gas pipeline networks. Various modeling approaches were introduced and explored for gas network flow, and were examined from analytical and numerical perspectives \cite{kazi2024modeling,srinivasan2022numerical,hari2021operation,dyachenko2017operator,tokareva2024stochastic, kiuchi1994implicit,herty2010new,zlotnik2015optimal}. The primary engineering challenge is to compensate for the decrease in pressure along the direction of flow that is caused by friction at the pipe's inner surface. The one-dimensional isothermal Euler equations, incorporating Darcy–Weisbach friction model, form a $2 \times 2$ system of coupled partial differential equations (PDEs) that represent the relevant conservation laws. These equations adequately describe gas flow dynamics along the pipe, accounting for gas momentum loss in the wave- and shock-free physical regime \cite{osiadacz1984simulation,misra2020monotonicity}.

The operation of natural gas transport networks presents complex challenges within the domain of control theory. The field of boundary stabilization of PDEs offers a rich theoretical framework to address these challenges \cite{colombo2009optimal,gugat2011gas,zhang2021boundary,auriol2020output,redaud2022stabilizing,anfinsen2019adaptive,yu2022traffic,gugat2020closed}. Lyapunov-based control \cite{gugat2011gas,gugat2020closed,zhang2021boundary}, optimal control \cite{colombo2009optimal,zlotnik2019optimal}, and PDE backstepping \cite{auriol2020output,redaud2022stabilizing,anfinsen2019adaptive,yu2022traffic} are among the control approaches that have been applied to gas flow control in pipelines using PDE models.  Methods that responsively compensate for unanticipated variations in boundary conditions caused by unscheduled changes in consumption, while maintaining pressure above minimum requirements, compel  ongoing research. 

In this paper, we examine the control of gas flowing through a pipeline subject to fluctuating outlet flow that we suppose is generated by a periodic linear dynamic system, which we use as a model of variable consumer demand. By appropriately selecting the parameters of the dynamic system, a desired periodic fluctuation profile can be achieved. First, we linearize the nonlinear isothermal Euler equations including Darcy-Weisbach friction term around an equilibrium solution and express the system in canonical form, resulting in a coupled $2 \times 2$ hyperbolic PDE. We then employ an observer-based PDE backstepping controller, adapted from \cite{aamo2012disturbance}, to regulate the outlet pressure to a setpoint by manipulating the inlet pressure, thus compensating for variations in outlet flow. Furthermore, we extend the observer-based controller to account for uncertain yet bounded fluctuations in the outlet flow. The derived controller manipulates the inlet pressure to regulate the outlet pressure to a neighborhood of the setpoint.

\section{Preliminaries}
Isothermal gas flow along a pipe, in a wave- and shock-free regime, can adequately be described by the one-dimensional isothermal Euler equations \cite{osiadacz1984simulation,misra2020monotonicity}: 
\begin{align}\label{aaq1}
    \rho_t+(\rho u)_x&=0,\\\label{qwopr}
    (\rho u)_t+(p+\rho u^2)_x&=-\frac{\lambda}{2D}\rho u\vert u\vert-\rho gh',
\end{align}
for all $x\in(0,\ell)$ and $t>0$, where
\begin{equation}
p(t,x)=ZRT\rho(t,x)=\sigma^2\rho(t,x),
\end{equation}
for all $x\in[0,\ell]$ and $t>0$,
with the boundary conditions
\begin{align}\label{aaz2}
\rho(t,0)&=U(t),\\\label{aa6}
    \rho(t,\ell)u(t,\ell)&=\varphi_{\text{L}}+s(t),
\end{align}
for all $t>0$. In equations \eqref{aaq1}-\eqref{aa6}, the variables \(\rho(t,x)\), \(u(t,x)\), and \(p(t,x)\) represent the density, velocity, and pressure profiles of the gas, respectively, and \(h(x)\) is the elevation of the pipe. The parameters involved are the internal pipe diameter \(D\), the pipe length \(\ell\), the gravitational acceleration \(g\), and the speed of sound \(\sigma=\sqrt{ZRT}\) in the gas, where \(Z\), \(R\), and \(T\) are the gas compressibility factor, specific gas constant, and absolute temperature, respectively.  We use the ideal gas equation of state, so that controlling density corresponds directly to controlling pressure. We consider \(U(t)\) in equation \eqref{aaz2} to be the control, which represents the density at the discharge of a compressor located at the pipe inlet.  The flow leaving the pipe outlet at \(x=\ell\) is \(\varphi_{\text{L}}>0\) in equation \eqref{aa6}, and \(s(t)\) in \eqref{aaz2} represents the variation in the outlet flow with respect to the nominal value \(\varphi_{\rm L}\).

\begin{ass}\label{ass1}
    The outlet flow variation $s(t)$ takes the form 
    \begin{equation}\label{st}
        s(t) = CX(t),
    \end{equation}
    where $X(t)\in\mathbb{R}^{n\times 1}$,
    \begin{equation} \label{dy_dis}
        \dot{X}(t) = AX(t),
    \end{equation}
 $C\in\mathbb{R}^{1\times n},$ and $A\in\mathbb{R}^{n\times n}$. The initial condition of $X(t)$ is in the set 
\begin{equation}\label{init_x}
    X(0) =\{X_0\vert CX_0=s(0)\}.
\end{equation}
\end{ass}
\begin{rmk}
Assumption \ref{ass1} is intended to model pipe outlet flow variation $s(t)$ to represent cyclic consumer energy demands. Periodic signals with specific amplitudes and periods can be constructed or approximated via \eqref{st}-\eqref{init_x} by choosing $A\in\mathbb{R}^{n\times n}$, $C$, and $n$, appropriately.
\end{rmk}
In the wave- and shock-free regime, the gas advection term $(\rho u^2)_x$ in \eqref{qwopr} is often ignored since it is very small when compared with the pressure gradient $p_x$ \cite{osiadacz1984simulation,misra2020monotonicity}. Further, we assume that the pipeline is level, which allows the gravity term $\rho gh'$ in \eqref{qwopr} to be ignored. Let us define 
\begin{equation}\label{zmlpdf}
    \varphi(t,x):=\rho(t,x)u(t,x),
\end{equation}
for all $x\in [0,\ell]$ and $t\geq 0$. Here, $\varphi(t,x)$ is the mass flux per area. This allows us to rewrite the system \eqref{aaq1}-\eqref{zmlpdf} as  
\begin{align}\label{qq1}
    \rho_t+\varphi_x&=0,\\
    \varphi_t+\sigma^2\rho_x&=-\frac{\lambda}{2D}\frac{\varphi\vert \varphi\vert}{\rho},
\end{align}
for all $x\in(0,\ell)$ and $t>0$ with the boundary conditions
\begin{align}\label{aazd2}
\rho(t,0)&=U(t),\\\label{aad6}\varphi(t,\ell)&=\varphi_{\text{L}}+s(t),
\end{align}
for all $t>0$.

\subsection{Equilibrium Steady-State Solution}
The system \eqref{qq1}-\eqref{aad6} with $U(t):=U_\star>0$ and $s(t)\equiv 0$ allows a continuum of equilibrium solutions, namely
\begin{align}
    \varphi(x)&\equiv \varphi_\star=\varphi_{\text{L}}>0,\\\label{xmlpr}
    \rho(x)&\equiv  \rho_\star(x)=\sqrt{(U_\star)^2-\frac{\lambda \varphi_{\text{L}}^2}{\sigma^2D} x}, \text{ }\forall x\in[0,\ell].  
\end{align}
As evident from equation \eqref{xmlpr}, the equilibrium control input $U_\star>0$ should be chosen such that 
\begin{align}\label{zmvlpdf}
    U_\star> \sqrt{\frac{\lambda \varphi_{\text{L}}^2\ell}{\sigma^2D} }.
\end{align}
\subsection{Linearization}
Let us define
\begin{align}
    \delta\rho(t,x)&:=\rho(t,x)-\rho_\star(x),\\
    \delta \varphi(t,x)&:=\varphi(t,x)-\varphi_{\text{L}},\\
    \delta U(t)&:=U(t)-U_\star,
\end{align}
for all $x\in[0,\ell]$ and $t\geq 0$, where $(\phi_{\rm L},\rho_{\star}(x),U_{\star})$ is the equilibrium. Further, because we consider small perturbations around the equilibrium $(\varphi_{\text{L}},\rho_\star(x),U_{\star})$ with $\varphi_{\text{L}}>0$, $\rho_\star(x)>0$ for all $x\in[0,\ell]$, and $U_{\star}>0$ chosen as in equation \eqref{zmvlpdf}, let us assume that $\varphi(t,x)>0$ for all $t> 0$ and $x\in [0,\ell]$. Then, linearizing the system \eqref{qq1}-\eqref{aad6} with $s(t)\equiv 0$ around the equilibrium $(\varphi_{\text{L}},\rho_\star(x),U_{\star})$, we obtain  
\begin{align}\label{zmlpw1er1}
    \delta \rho_t+\delta\varphi_x&=0,\\
\delta\varphi_t+\sigma^2\delta\rho_x&=\lambda_1(x)\delta \rho-\lambda_2(x)\delta\varphi,  
\end{align}
for all $x\in(0,\ell),t>0$, where
\begin{equation}
        \lambda_1(x):=\frac{\lambda \varphi_{\text{L}}^2}{2D\big(\rho_\star(x)\big)^2}\text{ and }\lambda_2(x):=\frac{\lambda \varphi_{\text{L}}}{D\rho_\star(x)},
\end{equation}
for all $x\in(0,\ell)$ with boundary conditions
\begin{align}
    \delta\rho(t,0)&=\delta U(t),\\\label{zmlpw1er2}
    \delta\varphi(t,\ell)&=0,
\end{align}
for all $t>0$.

\subsection{Transforming the System \eqref{zmlpw1er1}-\eqref{zmlpw1er2} into Canonical Form}
Let us consider the following change of coordinates: 
\begin{equation}\label{xbar}
    \bar{x}=\frac{\ell-x}{\ell},
\end{equation}
and the change of variables
\begin{align}\label{zzmlper}
  \!\! \!\! &\!v\Big(t,\frac{\ell-x}{\ell}\Big)\nonumber\\&=\!\sqrt{\!\frac{\rho_\star(x)}{\rho_\star(\ell)}}e^{\tfrac{\sigma}{\varphi_{\rm L}}\big(\rho_\star(x)-\rho_\star(\ell)\big)}\!\Big(\!\frac{1}{2}\delta\rho(t,x)\!-\!\frac{1}{2\sigma}\delta\varphi(t,x)\!\Big),\!\\\label{ddgh}\!\!&\!\!w\Big(t,\frac{\ell-x}{\ell}\Big) \nonumber\\
  \!\!\!\!&=\!\sqrt{\!\frac{\rho_\star(x)}{\rho_\star(\ell)}}e^{\!-\frac{\sigma}{\varphi_{\rm L}}\!\big(\rho_\star(x)-\rho_\star(\ell)\big)}\!\!\Big(\!\frac{1}{2}\delta\rho(t,x)\!+\!\frac{1}{2\sigma}\delta\varphi(t,x)\!\Big),\!
\end{align}
 for all $x\in[0,\ell]$ and $t\geq 0$. One can show that $(v,w)$ satisfy
\begin{align}\label{zmllw1}
    v_t+ \frac{\sigma}{\ell}v_{\bar{x}}&=-\mu_1(\bar{x})w,\\\label{zmllw22}
    w_t- \frac{\sigma}{\ell}w_{\bar{x}}&=\mu_2(\bar{x})v(t,\bar{x}),
\end{align}
for all $\bar{x}\in (0,1)$ and $t>0$, with the boundary conditions
\begin{align}\label{zzmlq1}
    v(t,0)&=w(t,0),\\\label{zmllw2}
    w(t,1)&= -r_1v(t,1)+r_2\delta U(t),
\end{align}
for all $t>0$. In equations \eqref{zmllw1} and \eqref{zmllw22}, $\mu_1(\bar{x})$ and $\mu_2(\bar{x})$ are given by 
\begin{align}
\begin{split}
&\mu_1(\bar{x})=\Bigg[\frac{\lambda\varphi^2_{\rm L}}{4\sigma D \Big(\rho_\star \big(\ell(1-\bar{x})\big)\Big)^2}\\\label{xxvbxcv2}&\quad-\frac{\lambda\varphi_{\rm L}}{2D\rho_\star\big(\ell(1-\bar{x})\big)}\Bigg]e^{\tfrac{2\sigma}{\varphi_{\rm L}}\Big(\rho_\star(\ell(1-\bar{x}))-\rho_\star(\ell)\Big)},\end{split}\\\begin{split}&\mu_2(\bar{x})=\Bigg[\frac{\lambda\varphi^2_{\rm L}}{4\sigma D \Big(\rho_\star \big(\ell(1-\bar{x})\big)\Big)^2}\\&\quad+\frac{\lambda\varphi_{\rm L}}{2D\rho_\star\big(\ell(1-\bar{x})\big)}\Bigg]e^{-\tfrac{2\sigma}{\varphi_{\rm L}}\Big(\rho_\star(\ell(1-\bar{x}))-\rho_\star(\ell)\Big)},
\end{split}
\end{align}
for all $\bar{x}\in(0,1)$. In equation \eqref{zmllw2}, $r_1$ and $r_2$ are given by
\begin{align}r_1&=e^{-\tfrac{2\sigma}{\varphi_{\text{L}}}\big(\rho_\star(0)-\rho_\star(\ell)\big)},\\\label{r2}
r_2&=\sqrt{\frac{\rho_\star(0)}{\rho_\star(\ell)}}e^{-\tfrac{\sigma}{\varphi_{\rm L}}\big(\rho_\star(0)-\rho_\star(\ell)\big)}.
\end{align}

The signs of the transport speed indicate that the variable $v$ represents information that travels from left to right and $w$ represents information that travel from right to left. The system \eqref{zmllw1}-\eqref{r2} is well-posed with boundary conditions on the left and right specified for $v$ and $w$, respectively. 

\begin{ass}\label{zzbmlpcv}
The outlet flow fluctuation $s(t)=CX(t)$ is such that $s(t)\ll \phi_L$.
\end{ass}

In light of Assumption \ref{zzbmlpcv}, when outlet flow fluctuation is present, the boundary condition \eqref{zmlpw1er2} is modified as
\begin{equation}
    \delta\varphi(t,\ell)=CX(t). 
\end{equation}
With this modification, we rewrite the linearized system in canonical form 
\begin{align}\label{zmllw1df}
    v_t+ \frac{\sigma}{\ell}v_{\bar{x}}&=-\mu_1(\bar{x})w,\\\label{zmllw22df}
    w_t- \frac{\sigma}{\ell}w_{\bar{x}}&=\mu_2(\bar{x})v,
\end{align}
for all $\bar{x}\in (0,1)$ and $t>0$, with the boundary conditions
\begin{align}\label{zzmlq1df}
    v(t,0)&=w(t,0)-\frac{1}{\sigma}CX(t),\\\label{zmllw2df}
    w(t,1)&= -r_1v(t,1)+r_2\delta U(t),
\end{align}
for all $t>0$, where $\mu_1(\bar{x})$, $\mu_2(\bar{x})$, $r_1$, and $r_2$ are as in equations \eqref{xxvbxcv2}-\eqref{r2}, respectively, and $X(t)$ satisfies eq. \eqref{dy_dis}. 

\subsection{Problem Formulation}
Our goal is to regulate the variation in outlet density, \(\delta \rho(t, \ell) \rightarrow 0\), despite the presence of outlet flow fluctuations, \(\delta \varphi(t, \ell) = CX(t)\). Note that by equations \eqref{zzmlper} and \eqref{ddgh},  
\begin{equation}
    \rho(t,\ell) = v(t,0)+w(t,0).
\end{equation}
Therefore, regulation of outlet density variation \(\delta \rho(t, \ell) \rightarrow 0\) is equivalent to \(v(t, 0) \rightarrow -w(t, 0)\). 
\section{Regulation of Outlet Density}

\subsection{When outlet flow varies according to Equations \eqref{st}-\eqref{init_x}}\label{sct23}
Similar to the approach in an earlier study \cite{aamo2012disturbance}, we use backstepping boundary control to achieve \(v(t,0) \rightarrow -w(t,0)\) subject to Assumption \ref{ass1} using the boundary measurement \(v(t,1)\). 

Let us consider the following observer
\begin{align}\label{dis_zmllw1sdf1}
    \hat{v}_t+ \frac{\sigma}{\ell}\hat{v}_{\bar{x}}&=-\mu_1(\bar{x})\hat{w}+p_1(\bar{x})\big(v(t,1)-\hat{v}(t,1)\big),\\
    \hat{w}_t- \frac{\sigma}{\ell}\hat{w}_{\bar{x}}&=\mu_2(\bar{x})\hat{v}+p_2(\bar{x})\big(v(t,1)-\hat{v}(t,1)\big),
\end{align}
for all $\bar{x}\in (0,1)$ and $t>0$,
\begin{align}\label{dis_zzmlq1f1}
    \hat{v}(t,0)&=\hat{w}(t,0)-\frac{1}{\sigma}CX(t),\\\label{zmllw2f1}
    \hat{w}(t,1)&= -r_1v(t,1)+r_2\delta U(t),
\end{align}
for all $t>0$, where $p_1(\bar{x})$ and $p_2(\bar{x})$ are observer gains chosen as
\begin{align}\label{ss11}
    p_1(\bar{x})&=-\frac{\sigma}{\ell}P^{11}(\bar{x},1),\\\label{ss21}
    p_2(\bar{x})&=-\frac{\sigma}{\ell}P^{21}(\bar{x},1),
\end{align}
for all $\bar{x}\in (0,1)$ with $P^{11}$ and $P^{21}$ satisfying 
\begin{align}\label{LL1ed}
 P^{11}_{\bar{x}}({\bar{x}},\xi) +  P^{11}_\xi({\bar{x}},\xi) &= -\frac{\ell}{\sigma}\mu_1(\bar{x}) P^{21}({\bar{x}},\xi), \\
 P^{21}_{\bar{x}}({\bar{x}},\xi)  -  P^{21}_\xi({\bar{x}},\xi) &=   -\frac{\ell}{\sigma}\mu_2(\bar{x}) P^{11}({\bar{x}},\xi),
\end{align}
with boundary conditions
\begin{align}
P^{11}(0,\xi) &=  P^{21}(0,\xi),\\\label{xxmlzpur}
P^{21}(\bar{x}, \bar{x}) &= -\frac{\ell}{2\sigma}\mu_2(\bar{x}).
\end{align}
The coupled hyperbolic PDEs \eqref{LL1ed}-\eqref{xxmlzpur} operate in the triangular domain $\{0\leq \bar{x} \leq \xi\leq 1\}$. Further, let us choose the control input $\delta U(t)$ as 
\begin{equation}\label{aazmlpsd}
\begin{split}
    \delta U(t) =& \frac{r_1}{r_2}v(t,1)+
    \frac{1}{r_2}\int_{0}^{1}K^{21}(1,\xi)\hat{v}(t,\xi)d\xi\\&+\frac{1}{r_2}\int_{0}^{1}K^{22}(1,\xi)\hat{w}(t,\xi)d\xi+\frac{1}{r_2}KX(t),
\end{split}
\end{equation}
for all $t>0$, where $K\in\mathbb{R}^{1\times n},K^{21},K^{22}$ are control gains satisfying 
\begin{equation}\label{kkmlpdd}
    K = \frac{1}{2\sigma}Ce^{A\frac{\ell}{\sigma}}-\frac{1}{\sigma}\int_{0}^{1} K^{21}(\tau,0)Ce^{A\frac{\ell}{\sigma}(1-\tau)}d\tau,
\end{equation}
and
\begin{align}\label{llljk1}
 K^{21}_{\bar{x}}({\bar{x}},\xi)  -  K^{21}_\xi({\bar{x}},\xi) &=   \frac{\ell}{\sigma}\mu_2(\xi) K^{22}({\bar{x}},\xi), \\\label{llljk1h}
 K^{22}_{\bar{x}}({\bar{x}},\xi)  +  K^{22}_\xi({\bar{x}},\xi) &=  - \frac{\ell}{\sigma}\mu_1(\xi) K^{21}({\bar{x}},\xi), 
\end{align}
with boundary conditions 
\begin{align}
K^{21}(\bar{x}, \bar{x}) &= -\frac{\ell}{2\sigma}\mu_2(\bar{x}), \\\label{kk4}
K^{22}(\bar{x}, 0) &=  K^{21}(\bar{x}, 0).
\end{align}
The coupled hyperbolic PDEs \eqref{llljk1}-\eqref{kk4} operate in the triangular domain $\{0\leq\xi\leq\bar{x}\leq 1\}$. Then, we can obtain the following result:

\begin{thme}[adapted from \cite{aamo2012disturbance}]
    Let the observer gains $p_1(\bar{x}),p_2(\bar{x})$ be chosen as in equations \eqref{ss11}-\eqref{xxmlzpur} and the control input $\delta U(t)$ be chosen as \eqref{aazmlpsd}-\eqref{kk4}. Then, $v(t,0)= -w(t,0)$ \textit{i.e} $\delta \rho(t,\ell)=0$ for all $t\geq 3\ell/\sigma$.  
\end{thme}

\subsection{When the outlet flow variation is subject to uncertainties}\label{sct33}

In the previous subsection, outlet density variation $\delta\rho(t,\ell)$ is regulated to zero in finite time when outlet flow is varying. The outlet flow variation is captured by $s(t)=CX(t)$ with the dynamics of $X(t)$ given by equation \eqref{dy_dis}. Here, we suppose that the flow variation $s(t)$ is subject to an unknown but bounded disturbance. That is, the outlet flow is given by
\begin{equation}
    \varphi(t,\ell) = \varphi_{\rm L}+d(t),
\end{equation}
where 
\begin{equation}
    d(t) = s(t) + \varepsilon(t),
\end{equation}
with $s(t)$ satisfying equations \eqref{st} and \eqref{dy_dis}, and $\varepsilon(t)$ is an unknown but bounded disturbance.
\begin{ass}\label{xxmlsd}
        The unknown disturbance $\varepsilon(t)\in C^{1}(\mathbb{R}_+)$ is bounded. That is, there exists a constant $M>0$ such that
    \begin{equation}
        \vert \varepsilon(t)\vert \leq M,
    \end{equation}
    for all $t\geq 0$. Furthermore, it holds that $M\ll\varphi_{\rm L}$. 
\end{ass}

In light of Assumption \ref{xxmlsd}, when outlet flow fluctuation with uncertainties is present, the boundary condition \eqref{zmlpw1er2} is modified as follows:
\begin{equation}
    \delta\varphi(t,\ell)=CX(t)+\varepsilon(t). 
\end{equation}

With this modification, we rewrite the linearized system in canonical form
\begin{align}\label{zmllw1dfdf}
    v_t+ \frac{\sigma}{\ell}v_{\bar{x}}&=-\mu_1(\bar{x})w,\\\label{zmllw22dfdf}
    w_t- \frac{\sigma}{\ell}w_{\bar{x}}&=\mu_2(\bar{x})v,
\end{align}
for all $\bar{x}\in (0,1)$ and $t>0$, with the boundary conditions
\begin{align}\label{zzmlq1dfdf}
    v(t,0)&=w(t,0)-\frac{1}{\sigma}CX(t)-\frac{1}{\sigma}\varepsilon(t),\\\label{zmllw2dfdf}
    w(t,1)&= -r_1v(t,1)+r_2\delta U(t),
\end{align}
for all $t>0$, where $\mu_1(\bar{x})$, $\mu_2(\bar{x})$, $r_1$, and $r_2$ are given by equations \eqref{xxvbxcv2}-\eqref{r2}, respectively, and $X(t)$ is generated by the system \eqref{dy_dis}. 

Let us consider the observer
\begin{align}\label{obst_dis_zmllw1sdf1}
    \hat{v}_t+ \frac{\sigma}{\ell}\hat{v}_{\bar{x}}&=-\mu_1(\bar{x})\hat{w}+p_1(\bar{x})\big(v(t,1)-\hat{v}(t,1)\big),\\\label{obst_dis_zmllw1sdf2}
    \hat{w}_t- \frac{\sigma}{\ell}\hat{w}_{\bar{x}}&=\mu_2(\bar{x})\hat{v}+p_2(\bar{x})\big(v(t,1)-\hat{v}(t,1)\big),
\end{align}
for all $\bar{x}\in (0,1)$ and $t>0$, with
\begin{align}\label{obst_dis_zzmlq1f1}
    \hat{v}(t,0)&=\hat{w}(t,0)-\frac{1}{\sigma}C\hat{X}(t),\\\label{obst_zmllw2f1}
    \hat{w}(t,1)&= -r_1v(t,1)+r_2\delta U(t),
\end{align}
for all $t>0$, and
\begin{align}\label{x_hat}
    \dot{\hat{X}}(t)=A\hat{X}(t)+ e^{A\frac{\ell}{\sigma}} H \big(v(t,1)-\hat{v}(t,1)\big),
\end{align}
for all $t>0$. Let the observer gains $H\in\mathbb{R}^{n\times 1},p_1(\bar{x}),$ and $p_2(\bar{x})$ be chosen such that $A+\frac{1}{\sigma}HC$ is hurwitz,
\begin{align}\label{zmldpuor1}
        p_1(\bar{x})&=-\frac{1}{\sigma}Ce^{A\frac{\ell}{\sigma}}H-\frac{\sigma}{\ell}P^{11}(\bar{x},1)\\&\nonumber+\frac{1}{\sigma}\int_{\bar{x}}^{1}P^{11}(\bar{x},\xi)Ce^{A \frac{\ell}{\sigma}(1-\xi)}Hd\xi,\\\label{zmldpuor2}
        p_2(\bar{x})&=-\frac{\sigma}{\ell}P^{21}(\bar{x},1)+\frac{1}{\sigma}\int_{\bar{x}}^{1}P^{21}(\bar{x},\xi)Ce^{A \frac{\ell}{\sigma}(1-\xi)}Hd\xi,
    \end{align}
where $P^{11}$ and $P^{21}$ are solutions to the system \eqref{LL1ed}-\eqref{xxmlzpur}. Further, let the control input $\delta U(t)$ be chosen as
\begin{equation}\label{aazmlpsdt1}
\begin{split}
    \delta U(t) =& \frac{r_1}{r_2}v(t,1)+
    \frac{1}{r_2}\int_{0}^{1}K^{21}(1,\xi)\hat{v}(t,\xi)d\xi\\&+\frac{1}{r_2}\int_{0}^{1}K^{22}(1,\xi)\hat{w}(t,\xi)d\xi+\frac{1}{r_2}K\hat{X}(t),
\end{split}
\end{equation}
where $K$ is given by eq. \eqref{kkmlpdd}, and $K^{21},K^{22}$ satisfy equations \eqref{llljk1}-\eqref{kk4}. Then, we can obtain the following result. 
\begin{thme}\label{thm66}
    Let the observer gain vector $H\in\mathbb{R}^{n\times 1}$ in \eqref{x_hat} is chosen such that $A+\frac{1}{\sigma}HC$ is hurwitz. Further, let the observer gain functions $p_1(\bar{x})$ and $p_2(\bar{x})$ in equations \eqref{obst_dis_zmllw1sdf1} and \eqref{obst_dis_zmllw1sdf2} be set as equations \eqref{zmldpuor1}-\eqref{zmldpuor2}, and let the control input $\delta U(t)$ be chosen as equation \eqref{aazmlpsdt1}. Then, $\delta\rho(t,\ell)$ converges to a neighborhood of the origin as $t\rightarrow \infty$, i.e.,
\begin{equation}
    \lim_{t\rightarrow\infty}\vert \delta\rho(t,\ell)\vert \leq \bar{M},
\end{equation}
for some $\bar{M}>0$.  
\end{thme}
\noindent \textbf{Proof.} The proof contains two parts: I) obtaining an expression for $\delta\rho(t,\ell)$ valid for all $t\geq \ell/\sigma$; and II) obtaining an upper-bound for the limit $\lim_{t\rightarrow\infty}\vert\delta\rho(t,\ell)\vert$. 

\noindent \underline{I) An expression for $\delta\rho(t,\ell)$ valid for all $t\geq \ell/\sigma$ }

Define
\begin{align}\label{zmlpwcdf1}
    \tilde{v}(t,\bar{x}) &:= v(t,\bar{x})-\hat{v}(t,\bar{x}),\\
    \tilde{w}(t,\bar{x}) &:= w(t,\bar{x})-\hat{w}(t,\bar{x}),\\\label{zmlpwcdf3}
    \tilde{X}(t)&:=X(t)-\hat{X}(t),
\end{align}
for all $\bar{x}\in [0,1]$ and $t\geq 0$. Then, we can rewrite the control input $\delta U(t)$ chosen in equation \eqref{aazmlpsdt1} as,  
\begin{align}
    \delta U(t) =&\frac{r_1}{r_2}v(t,1)+
    \frac{1}{r_2}\int_{0}^{1}K^{21}(1,\xi)v(t,\xi)d\xi\\&+\frac{1}{r_2}\int_{0}^{1}K^{22}(1,\xi)w(t,\xi)d\xi+\frac{1}{r_2}V(t),
\end{align}
where
\begin{equation}\label{zmcnvlpirt}
\begin{split}
    V(t):=&KX(t)-\int_{0}^{1}K^{21}(1,\xi)\tilde{v}(t,\xi)d\xi\\&-\int_{0}^{1}K^{22}(1,\xi)\tilde{w}(t,\xi)d\xi-K\tilde{X}(t).
\end{split}
\end{equation}
Consider the following backstepping transformations:

\begin{align}\label{bt1sd}
\begin{split}
    \alpha(t,\bar{x}) &= v(t,\bar{x})-\int_{0}^{\bar{x}}K^{11}(\bar{x},\xi)v(t,\xi)d\xi\\&\quad-\int_{0}^{\bar{x}}K^{12}(\bar{x},\xi)w(t,\xi)d\xi,
    \end{split}
    \\\label{bt2sd}
    \begin{split}
     \beta(t,\bar{x}) &= w(t,\bar{x})-\int_{0}^{\bar{x}}K^{21}(\bar{x},\xi)v(t,\xi)d\xi\\&\quad-\int_{0}^{\bar{x}}K^{22}(\bar{x},\xi)w(t,\xi)d\xi,
    \end{split}
\end{align}
defined in the triangular domain $\{0\leq \xi\leq \bar{x}\leq 1\}$, where $K^{21}$ and $K^{22}$ satisfy equations \eqref{llljk1}-\eqref{kk4}, and $K^{11}$ and $K^{12}$ satisfy
\begin{align}\label{lll1cxz}
 K^{11}_{\bar{x}}({\bar{x}},\xi) +  K^{11}_\xi({\bar{x}},\xi) &= - \frac{\ell}{\sigma}\mu_2(\xi) K^{12}({\bar{x}},\xi),  \\
 K^{12}_{\bar{x}}({\bar{x}},\xi)  -  K^{12}_\xi({\bar{x}},\xi) &=   \frac{\ell}{\sigma}\mu_1(\xi) K^{11}({\bar{x}},\xi), 
\end{align}
with boundary conditions 
\begin{align}
K^{11}(\bar{x}, 0) &=  K^{12}(\bar{x}, 0), \\
K^{12}(\bar{x}, \bar{x}) &= -\frac{\ell}{2\sigma}\mu_1(\bar{x}),
\end{align}
in the domain $\{0\leq \xi\leq \bar{x}\leq 1\}$.  Subject to the backstepping transformations \eqref{bt1sd} and \eqref{bt2sd}, the system \eqref{zmllw1dfdf}-\eqref{zmllw2dfdf} is transformed to the system
\begin{align}
    \alpha_t+\frac{\sigma}{\ell}\alpha_{\bar{x}}&=\frac{1}{\ell}K^{11}(\bar{x},0)CX(t)+\frac{1}{\ell}K^{11}(\bar{x},0)\varepsilon(t),\\
    \beta_t-\frac{\sigma}{\ell}\beta_{\bar{x}}&=\frac{1}{\ell}K^{21}(\bar{x},0)CX(t)+\frac{1}{\ell}K^{21}(\bar{x},0)\varepsilon(t),
\end{align}
for all $\bar{x}\in(0,1)$ and $t>0$, and the boundary conditions
\begin{align}
    \alpha(t,0)&=\beta(t,0)-\frac{1}{\sigma}CX(t)-\frac{1}{\sigma}\varepsilon(t),\\
    \beta(t,1)&=V(t),
\end{align}
for all $t>0$. Using the same line of reasoning as in proof of Lemma 3 and Theorem 4 of \cite{aamo2012disturbance}, we can show that 
\begin{align}\label{sscv1d}
\begin{split}
    &\alpha(t,0)=-\beta(t,0)+2V\Big(t-\frac{\ell}{\sigma}\Big)\\&\qquad+\frac{2}{\sigma}\int_{0}^{1}K^{21}(\tau,0)Ce^{-A\frac{\ell}{\sigma}\tau}d\tau X(t)-\frac{C}{\sigma}X(t)\\&\qquad
    +\frac{2}{\sigma}\int_{0}^{1}K^{21}(\tau,0)\varepsilon\Big(t-\frac{\ell}{\sigma}\tau\Big)d\tau-\frac{1}{\sigma}\varepsilon(t),
\end{split}
\end{align}
for all $t\geq \ell/\sigma$. Considering \eqref{zmcnvlpirt}, we can rewrite \eqref{sscv1d} as 
\begin{align}\label{zmkfptuo}
    \begin{split}
        &\alpha(t,0) = -\beta(t,0)+2KX\Big(t-\frac{\ell}{\sigma}\Big)\\&-2\int_{0}^{1}K^{21}(1,\xi)\tilde{v}\Big(t-\frac{\ell}{\sigma}\Big)d\xi\\&-2\int_{0}^{1}K^{22}(1,\xi)\tilde{w}\Big(t-\frac{\ell}{\sigma},\xi\Big)d\xi-2K\tilde{X}\Big(t-\frac{\ell}{\sigma}\Big)\\&+\frac{2}{\sigma}\int_{0}^{1}K^{21}(\tau,0)Ce^{-A\frac{\ell}{\sigma}\tau}d\tau X(t)-\frac{C}{\sigma}X(t)
    \\&+\frac{2}{\sigma}\int_{0}^{1}K^{21}(\tau,0)\varepsilon\Big(t-\frac{\ell}{\sigma}\tau\Big)d\tau-\frac{1}{\sigma}\varepsilon(t). 
    \end{split}
\end{align}
However, by the semigroup property of the system \eqref{dy_dis}, we have that $X(t-\ell/\sigma)=e^{-A\frac{\ell}{\sigma}}X(t)$. Using this fact and rearranging the terms of equation \eqref{zmkfptuo}, we can obtain
\begin{align}\label{nnmlpklpij}
\begin{split}
    &\alpha(t,0) =-\beta(t,0)+2Ke^{-A\frac{\ell}{\sigma}}X(t)\\&+\frac{2}{\sigma}\int_{0}^{1}K^{21}(\tau,0)Ce^{-A\frac{\ell}{\sigma}\tau}d\tau X(t)-\frac{C}{\sigma}X(t)\\&-2\int_{0}^{1}K^{21}(1,\xi)\tilde{v}\Big(t-\frac{\ell}{\sigma}\Big)d\xi\\&-2\int_{0}^{1}K^{22}(1,\xi)\tilde{w}\Big(t-\frac{\ell}{\sigma},\xi\Big)d\xi-2K\tilde{X}\Big(t-\frac{\ell}{\sigma}\Big)
    \\&+\frac{2}{\sigma}\int_{0}^{1}K^{21}(\tau,0)\varepsilon\Big(t-\frac{\ell}{\sigma}\tau\Big)d\tau-\frac{1}{\sigma}\varepsilon(t),
\end{split}
\end{align}
for all $t\geq \ell/\sigma$. Recalling that $K$ is chosen as in equation \eqref{kkmlpdd}, we can simplify \eqref{nnmlpklpij} to obtain that
\begin{align}\label{nnmlp}
\begin{split}
    &\alpha(t,0) = -\beta(t,0)-2\int_{0}^{1}K^{21}(1,\xi)\tilde{v}\Big(t-\frac{\ell}{\sigma}\Big)d\xi\\&-2\int_{0}^{1}K^{22}(1,\xi)\tilde{w}\Big(t-\frac{\ell}{\sigma},\xi\Big)d\xi-2K\tilde{X}\Big(t-\frac{\ell}{\sigma}\Big)
    \\&+\frac{2}{\sigma}\int_{0}^{1}K^{21}(\tau,0)\varepsilon\Big(t-\frac{\ell}{\sigma}\tau\Big)d\tau-\frac{1}{\sigma}\varepsilon(t),
\end{split}
\end{align}
for all $t\geq \ell/\sigma$. It follows from  equations \eqref{zzmlper} and \eqref{ddgh} that
\begin{align}
    \delta \rho(t,\ell) = v(t,0)+w(t,0).
\end{align}
However, equations \eqref{bt1sd} and \eqref{bt2sd} imply that $\alpha(t,0)=v(t,0)$ and $\beta(t,0)=w(t,0)$, and thus equation \eqref{nnmlp} yields
\begin{align}\label{nnmslpdf}
\begin{split}
    &\delta \rho(t,\ell) = -2\int_{0}^{1}K^{21}(1,\xi)\tilde{v}\Big(t-\frac{\ell}{\sigma},\xi\Big)d\xi\\&-2\int_{0}^{1}K^{22}(1,\xi)\tilde{w}\Big(t-\frac{\ell}{\sigma},\xi\Big)d\xi-2K\tilde{X}\Big(t-\frac{\ell}{\sigma}\Big)
    \\&+\frac{2}{\sigma}\int_{0}^{1}K^{21}(\tau,0)\varepsilon\Big(t-\frac{\ell}{\sigma}\tau\Big)d\tau-\frac{1}{\sigma}\varepsilon(t),
\end{split}
\end{align}
for all $t\geq \ell/\sigma$.

In order to analyse $\lim_{t\rightarrow\infty}\vert\delta\rho(t,\ell)\vert$, we examine $\lim_{t\rightarrow \infty}\tilde{v}(t)$, $\lim_{t\rightarrow\infty} \tilde{w}(t)$, and $\lim_{t\rightarrow\infty}\tilde{X}(t)$ as follows.

\noindent \underline{II) An upper-bound for $\lim_{t\rightarrow\infty}\vert \delta\rho(t,\ell)\vert$}

Considering equations \eqref{dy_dis}, \eqref{zmllw1dfdf}-\eqref{x_hat}, and \eqref{zmlpwcdf1}-\eqref{zmlpwcdf3}, we can show that the observer errors $(\tilde{v},\tilde{w},\tilde{X})$ satisfy   
\begin{align}\label{est_zmllw1sdf}
\tilde{v}_t+\frac{\sigma}{\ell} \tilde{v}_{\bar{x}}&=-\mu_1(\bar{x})\tilde{w}-p_1(\bar{x})\tilde{v}(1,t),\\
    \tilde{w}_t- \frac{\sigma}{\ell}\tilde{w}_{\bar{x}}&=\mu_2(\bar{x})\tilde{v}-p_2(\bar{x})\tilde{v}(1,t),
\end{align}
for all $\bar{x}\in (0,1)$ and $t>0$, with the boundary conditions
\begin{align}\label{est_zzmlq1f}
    \tilde{v}(t,0)&=\tilde{w}(t,0)-\frac{1}{\sigma}C\tilde{X}(t)-\frac{1}{\sigma}\varepsilon(t),\\\label{est_zmllw2f}
    \tilde{w}(t,1)&=0,
\end{align}
for all $t>0$, and 
\begin{equation}
    \dot{\tilde{X}}(t) = A\tilde{X}(t)-e^{A\frac{\ell}{\sigma}}H\tilde{v}(t,1),
\end{equation}
for all $t>0$. Consider the backstepping transformations
\begin{align}\label{wwel}
\begin{split}
    \tilde{v}(t,\bar{x}) &= \tilde{\alpha}(t,\bar{x})-\int_{\bar{x}}^{1}P^{11}(\bar{x},\xi)\tilde{\alpha}(t,\xi)d\xi\\&\quad-\int_{\bar{x}}^{1}P^{12}(\bar{x},\xi)\tilde{\beta}(t,\xi)d\xi,
\end{split}\\\label{wwel2}
\begin{split}
     \tilde{w}(t,\bar{x}) &= \tilde{\beta}(t,\bar{x})-\int_{\bar{x}}^{1}P^{21}(\bar{x},\xi)\tilde{\alpha}(t,\xi)d\xi\\&\quad-\int_{\bar{x}}^{1}P^{22}(\bar{x},\xi)\tilde{\beta}(t,\xi)d\xi,
\end{split}
\end{align}
defined in the triangular domain $\{0\leq \bar{x} \leq \xi\leq 1\}$, where $P^{11}$ and $P^{21}$ are solutions to equations \eqref{LL1ed}-\eqref{xxmlzpur}, and $P^{21}$ and $P^{22}$ satisfy 
\begin{align}\label{LL1edgh}
 P^{12}_{\bar{x}}({\bar{x}},\xi)  -  P^{12}_\xi({\bar{x}},\xi) &=  -\frac{\ell}{\sigma}\mu_1(\bar{x}) P^{22}({\bar{x}},\xi),  \\
 P^{22}_{\bar{x}}({\bar{x}},\xi)  +  P^{22}_\xi({\bar{x}},\xi) &=  - \frac{\ell}{\sigma}\mu_2(\bar{x}) P^{12}({\bar{x}},\xi), 
\end{align}
with boundary conditions
\begin{align}
P^{12}(\bar{x}, \bar{x}) &= -\frac{\ell}{2\sigma}\mu_1(\bar{x}), \\\label{LL55ed}
P^{22}(0,\xi) &=  P^{12}(0,\xi).
\end{align}
Subject to the observer error backstepping transformations \eqref{wwel} and \eqref{wwel2}, and the observer gain functions $p_1(\bar{x}),p_2(\bar{x})$ chosen as in equations \eqref{zmldpuor1} and \eqref{zmldpuor2}, the observer error system is transformed to 
\begin{align}
\tilde{\alpha}_t(t,\bar{x})+\frac{\sigma}{\ell}\tilde{\alpha}_{\bar{x}}(t,\bar{x})&=\frac{1}{\sigma}Ce^{A\frac{\ell}{\sigma}(1-\bar{x})}H\tilde{\alpha}(t,1),\\
    \tilde{\beta}(t,\bar{x})-\frac{\sigma}{\ell}\tilde{\beta}_{\bar{x}}(t,\bar{x})&=0,
\end{align}
for all $\bar{x}\in (0,1)$ and $t>0$ and the boundary conditions 
\begin{align}
     \tilde{\alpha}(t,0)&=\tilde{\beta}(t,0)-\frac{1}{\sigma}C\tilde{X}(t)-\frac{1}{\sigma}\varepsilon(t),\\
    \tilde{\beta}(t,1)&=0,
\end{align}
for all $t>0$, and 
\begin{align}
    \dot{\tilde{X}}(t)=A\tilde{X}(t)-e^{A\tfrac{\ell}{\sigma}}H\tilde{\alpha}(t,1),
\end{align}
for all $t>0$. Notice that $\tilde{\beta}[t]=0$ for all $t\geq \ell/\sigma$. Therefore, it holds that
\begin{align}
\tilde{\alpha}_t(t,\bar{x})+\frac{\sigma}{\ell}\tilde{\alpha}_{\bar{x}}(t,\bar{x})&=\frac{1}{\sigma}Ce^{A\tfrac{\ell}{\sigma}(1-\bar{x})}H\tilde{\alpha}(t,1),\\
    \tilde{\alpha}(t,0)&=-\frac{1}{\sigma}C\tilde{X}(t)-\frac{1}{\sigma}\varepsilon(t),
\end{align}
for all $\bar{x}\in (0,1)$ and for all $t\geq \ell/\sigma$, and
\begin{align}
    \dot{\tilde{X}}(t)=A\tilde{X}(t)-e^{A\tfrac{\ell}{\sigma}}H\tilde{\alpha}(t,1),
\end{align}
for all $t>0$. Define
\begin{equation}\label{aamlprt}
    \gamma(t,\bar{x}):=\tilde{\alpha}(t,\bar{x})+\frac{1}{\sigma}Ce^{-A\tfrac{\ell}{\sigma}\bar{x}}\tilde{X}(t), 
\end{equation}
for all $\bar{x}\in (0,1)$ and for all $t\geq \ell/\sigma$.  Then, we can show that
\begin{align}
    \gamma_t(t,\bar{x})+\frac{\sigma}{\ell}\gamma_{\bar{x}}(t,\bar{x})&=0,
\end{align}
for all $\bar{x}\in(0,1)$ and $t\geq \ell/\sigma$ and the boundary conditions
\begin{align}
     \gamma(t,0)&=-\frac{1}{\sigma}\varepsilon(t),
\end{align}
for all $t\geq \ell/\sigma$. Further, we have that
\begin{align}
\begin{split}
    \dot{\tilde{X}}(t)=& \Big(A+\frac{1}{\sigma}e^{A\frac{\ell}{\sigma}}HCe^{-A\tfrac{\ell}{\sigma}}\Big)\tilde{X}(t)-e^{A\tfrac{\ell}{\sigma}}H\gamma(t,1)\\=& e^{A\frac{\ell}{\sigma}}\Big(A+\frac{1}{\sigma}HC\Big)e^{-A\tfrac{\ell}{\sigma}}\tilde{X}(t)-e^{A\tfrac{\ell}{\sigma}}H\gamma(t,1),
\end{split}
\end{align}
for all $t\geq \ell/\sigma$. 
Let 
\begin{align}
    \tilde{A} = e^{A\tfrac{\ell}{\sigma}}\Big(A+\frac{1}{\sigma}HC\Big)e^{-A\tfrac{\ell}{\sigma}}.
\end{align}
The gain vector $H$ is chosen such that $A+\frac{1}{\sigma}HC$ is hurwitz. Therefore, we have that $\tilde{A}$ also hurwitz. 

Referring to Proposition  3.2 of \cite{karafyllis2019input}, we can obtain the input-to-state stability results
\begin{equation}
   \!\!\! \Vert \gamma[t]\Vert \!\!\leq  e^{-\Big(t-\tfrac{5\ell}{2\sigma}\Big)} \Vert \gamma[\ell/\sigma]\Vert\!\!+\!\!\sqrt{\!\frac{1}{2\sigma\ell}}e^{\tfrac{3\ell}{2\sigma}}\!\!\max_{\ell/\sigma\leq \tau\leq t}\vert \varepsilon(\tau)\vert,
   \end{equation}
   \begin{equation}
\vert \gamma(t,1)\vert \leq  e^{-\Big(t-\tfrac{3\ell}{2\sigma}\Big)}\vert \max_{0\leq \bar{x}\leq 1}\vert \gamma(\ell/\sigma,\bar{x})\vert \!\!+\!\!\frac{e^{\tfrac{2\ell}{\sigma}}}{\sigma}\!\!\max_{\ell/\sigma \leq \tau\leq t}\vert \varepsilon(\tau)\vert,
\end{equation}
for all $t\geq \ell/\sigma$, and that
\begin{equation}
\begin{split}
        \Vert \tilde{X}(t)\Vert \leq & \Omega e^{-\upsilon \big(t-\frac{\ell}{\sigma}\big)}\Vert \tilde{X}\big(\ell/\sigma\big)\Vert \\&+\frac{\Omega}{\upsilon} e^{\Vert A\Vert \frac{\ell}{\sigma}}\Vert H\Vert \max_{\ell/\sigma\leq \tau\leq t}\vert \gamma (\tau,1)\vert, 
\end{split}
\end{equation}
where $\upsilon$ is the smallest absolute value of the real parts of the eigenvalues of the Hurwitz matrix $A+\frac{1}{\sigma}HC$. Therefore, as $t\rightarrow \infty$, we have that 
\begin{equation}\label{zzc1}
\!\lim_{t\rightarrow\infty}\Vert \gamma[t]\Vert \!=\! \sqrt{\!\frac{1}{2\sigma\ell}}e^{\tfrac{3\ell}{2\sigma}}\!\!\max_{\ell/\sigma\leq \tau\leq \infty}\vert \varepsilon(\tau)\vert \leq \sqrt{\frac{1}{2\sigma\ell}}e^{\tfrac{3\ell}{2\sigma}} M,
\end{equation}
\begin{align}\label{zzc2}\lim_{t\rightarrow\infty} \vert \gamma(t,1)\vert &= e^{\tfrac{2\ell}{\sigma}}\max_{\ell/\sigma \leq \tau\leq \infty}\vert \varepsilon(\tau)\vert/\sigma \leq e^{\tfrac{2\ell}{\sigma}} M/\sigma,\\\label{zzc3}
  \lim_{t\rightarrow\infty} \Vert \tilde{X}(t)\Vert &= \frac{\Omega}{\upsilon} e^{\Vert A\Vert\tfrac{\ell}{\sigma}}\Vert H\Vert \max_{\ell/\sigma\leq \tau\leq \infty}\vert\gamma (\tau,1)\vert\\&\leq \frac{\Omega}{\upsilon\sigma} e^{\Vert A\Vert\tfrac{\ell}{\sigma}}\Vert H\Vert e^{\tfrac{2\ell}{\sigma}} M. 
\end{align}
From equation \eqref{aamlprt}, we can obtain that
\begin{align}
\begin{split}
    \Vert \tilde{\alpha}[t]\Vert &\leq \Vert \gamma[t]\Vert + \frac{1}{\sigma}\Vert C\Vert e^{\Vert A\Vert\frac{\ell}{\sigma}} \Vert \tilde{X}(t)\Vert,
\end{split}
\end{align}
for all $t\geq \ell/\sigma$. Therefore, considering equations \eqref{zzc1} and \eqref{zzc3}, we can show that
\begin{equation}\label{rrrtyu}
  \tilde{\alpha}[t]\Vert \leq \sqrt{\frac{1}{2\sigma\ell}}e^{\tfrac{3\ell}{2\sigma}} M+\frac{\Omega}{\sigma^2\upsilon}\Vert C\Vert e^{\Vert A\Vert\frac{2\ell}{\sigma}} \Vert H\Vert e^{\tfrac{2\ell}{\sigma}} M,
\end{equation}
as $t\rightarrow\infty$. Recall that $\tilde{\beta}[t]=0$ for all $t\geq \ell/\sigma$. Therefore, considering the observer error backstepping transformations \eqref{wwel} and \eqref{wwel2}, we can obtain that 
\begin{equation}
    \Vert \tilde{v}[t] \Vert \leq \tilde{P}_{11}\Vert\tilde{\alpha}[t]\Vert,  \quad \text{ and } \quad
    \Vert \tilde{w}[t]\Vert \leq \tilde{P}_{21}\Vert\tilde{\alpha}[t]\Vert,  \label{pplklo}
\end{equation}
for all $t\geq \ell/\sigma$, where 
\begin{align}\label{tildeldfd}
\tilde{P}_{11}=&1+\Big(\int_{0}^{1}\int_{\bar{x}}^1P_{11}^2(\bar{x},\xi)d\xi d\bar{x}\Big)^{1/2},\\
\tilde{P}_{21}=& \int_{0}^{1}\int_{\bar{x}}^1P_{21}^2(\bar{x},\xi)d\xi d\bar{x}. 
\end{align}
Further, considering equation \eqref{nnmslpdf}, we can obtain
\begin{align}\label{bbbnk}
    \begin{split}
           & \vert \delta \rho(t,\ell)\vert  \leq  2\max_{0\leq \xi \leq 1}\vert K^{21}(1,\xi)\vert\cdot \Vert \tilde{v}[t-\ell/\sigma]\Vert \\&+2\max_{0\leq \xi\leq 1}\vert K^{22}(1,\xi)\vert \cdot \Vert\tilde{w}[t-\ell/\sigma]\Vert+\frac{M}{\sigma}\\&+2\Vert K\Vert\cdot\Vert\tilde{X}[t-\ell/\sigma]\Vert+\frac{2}{\sigma} \max_{(0\leq \xi\leq 1)}\vert K^{21}(\xi,0)\vert M , 
    \end{split}
\end{align}
for all $t\geq \ell/\sigma$. 
Then, using equations \eqref{zzc1}-\eqref{zzc3}, \eqref{rrrtyu}-\eqref{pplklo}, and \eqref{bbbnk}, we can show that $\vert\delta\rho(t,\ell)\vert\rightarrow \bar{M}$ as $t\rightarrow \infty$, for some $\bar{M}>0$. This completes the proof. \hfill $\square$

\section{Numerical Simulations}
\begin{figure}[h!]
    \centering
    \includegraphics[width=\linewidth, height=4cm]{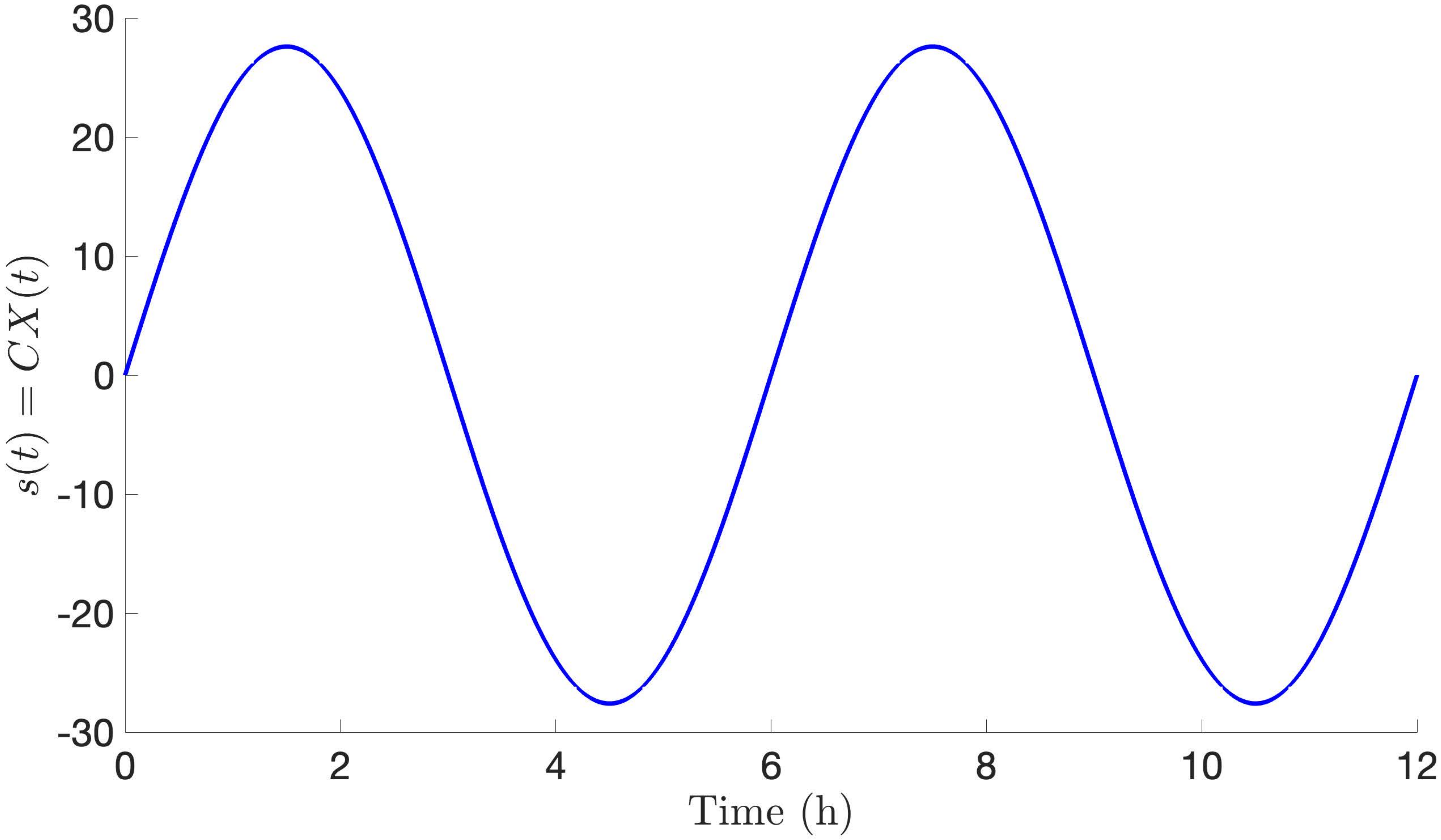}
    \caption{Variation in pipe outlet flow.}
    \label{fg1}
    \vspace{-3ex}
\end{figure}

We demonstrate the observer-based boundary control method developed above for setpoint tracking and disturbance attenuation  in a simple gas pipeline model.  Motivated by previous studies \cite{zlotnik2015model,tokareva2024stochastic}, we use the plant parameters $\lambda = 0.011$, $\sigma = 378$ m/s, $D = 0.5$ m, $\ell = 25\times 10^3$ m, $\varphi_{\rm L} = 289$ kg/m$^2$/s, and $U_\star = 46$ kg/m$^3$.
Let $A$ and $C$ be
\begin{align}
    A &= \begin{pmatrix}
        0 && 1\\-\big(\frac{2\pi}{6\times 3600}\big)^2 &&  0
    \end{pmatrix}, \text{ }
   C&=\begin{pmatrix}
        1&&0
    \end{pmatrix}. 
\end{align}
Further, let the initial conditions be chosen as  
\begin{align}\label{eq_rho_dis}
\rho(0,x) &= \rho_\star(x),\\
\phi(0,x) &= \phi_{\rm L},\\
X(0) & = \begin{pmatrix}
    0\\\frac{0.1\phi_{\rm L}}{3600}
\end{pmatrix}.
\end{align}
Then, considering the solution of the system \eqref{dy_dis} and \eqref{init_x}, we can show that $s(t)$ is given by
\begin{align}
    s(t) = \frac{0.6\varphi_{\rm L}}{2\pi}\sin\Big(\frac{2\pi t}{6\times 3600}\Big).
\end{align}
Fig. \ref{fg1} shows the evolution of $s(t)$. The initial conditions for the observer \eqref{dis_zmllw1sdf1}-\eqref{zmllw2f1}  are chosen as 
\begin{equation}
    \hat{v}[0]\equiv 0, \text{ }\hat{w}[0]\equiv 0. 
\end{equation}
Below, we conduct numerical simulations for two cases: A) when the outlet flow varies according to equations \eqref{st}-\eqref{init_x} as discussed in Section \ref{sct23}; and B) when the outlet flow variation is subject to uncertainties as in Section \ref{sct33}.

\subsection{When outlet flow varies according to \eqref{st}-\eqref{init_x}}
\begin{figure}[h!]
    \centering
    \includegraphics[width=\linewidth]{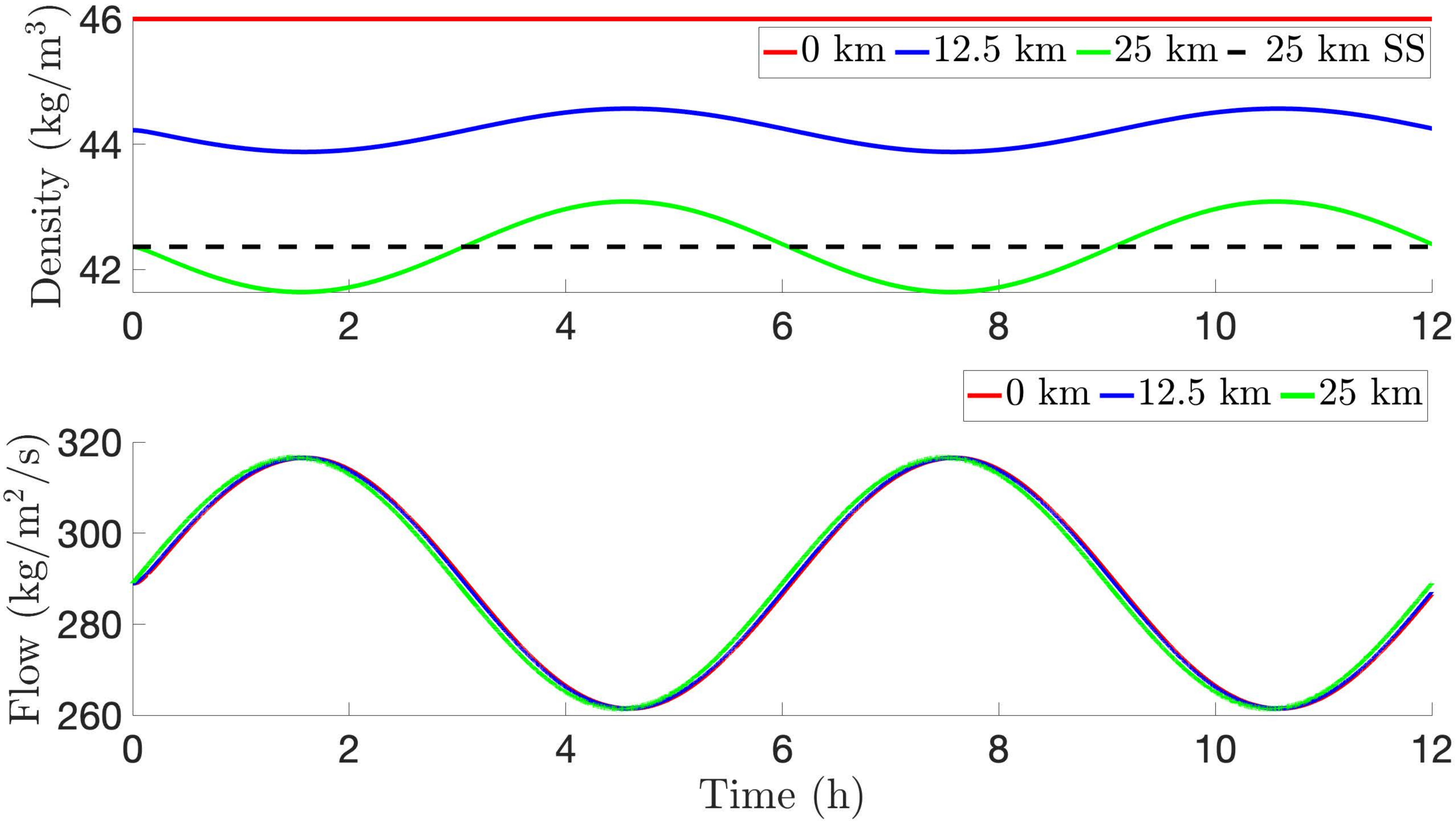}
    \caption{Density (top) and flow (bottom) with $\delta U(t)=0$}
    \label{fg2}
\end{figure}

\begin{figure}[h!]
    \centering
    \includegraphics[width=\linewidth]{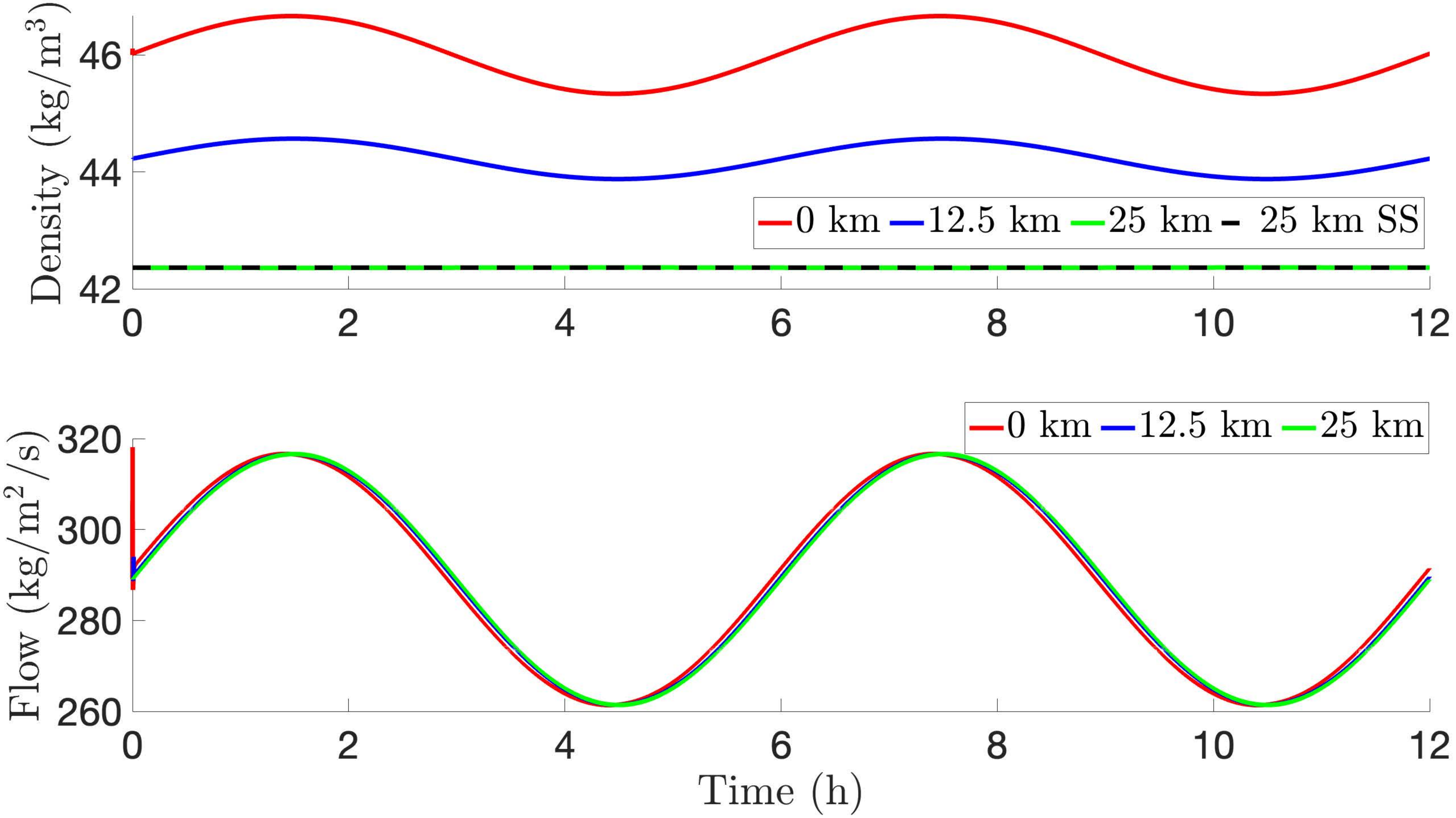}
    \caption{Density (top) and flow (bottom) with $\delta U(t)$ in \eqref{aazmlpsd}.}
    \label{fg3}
\end{figure}

In Figs. \ref{fg2} and \ref{fg3}, we illustrate the evolution of density and flow at three locations along the pipe. Fig. \ref{fg2} depicts the scenario where $\delta U(t) = 0$, while Fig. \ref{fg3} shows the case where $\delta U(t)$ is given by equation \eqref{aazmlpsd}. As seen in Fig. \ref{fg3}, the outlet density is regulated to the steady-state value even in the presence of outlet flow fluctuations. In contrast, Fig. \ref{fg2} shows that the outlet density fluctuates due to the presence of outlet flow disturbances.


\begin{figure}[h!]
    \centering
    \includegraphics[width=\linewidth]{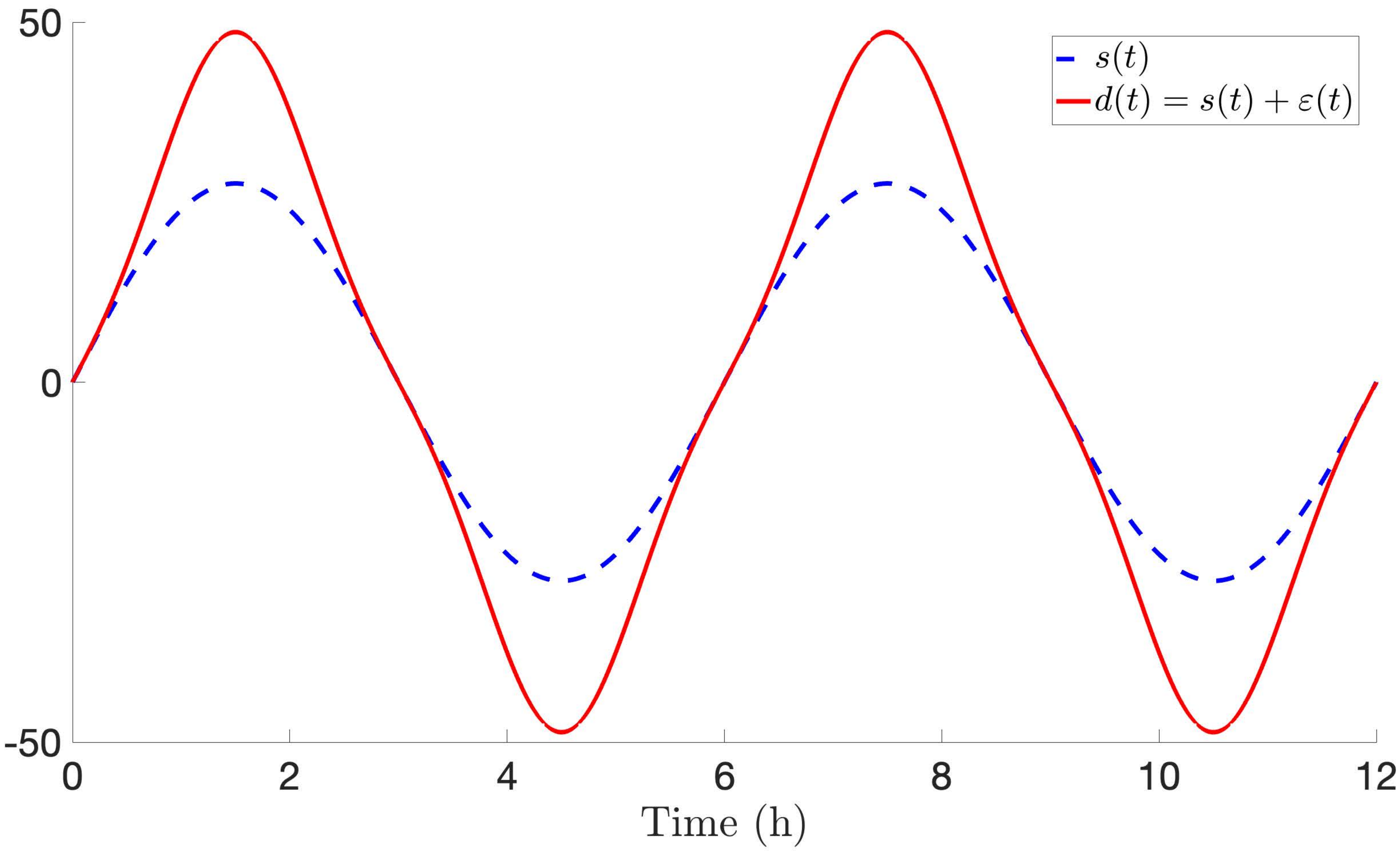}
    \caption{Variation in outlet flow with disturbance.}
    \label{fg5}
\end{figure}

\subsection{When the outlet flow variation is subject to uncertainties}

Here, we assume that the outlet flow variation is subject to an unknown yet bounded disturbance, $\varepsilon(t)$. For simulation purposes, we take
\begin{equation}
\varepsilon(t) = 0.001s^3(t).
\end{equation}
In Fig. \ref{fg5}, we show the evolution of the disturbed outlet flow. The initial conditions of the observer $\hat{X}(t)$ are chosen as
\begin{align}
\hat{X}(0) = \begin{pmatrix}
0 && 0
\end{pmatrix}^{T}.
\end{align}

\begin{figure}[h!]
    \centering
    \includegraphics[width=\linewidth]{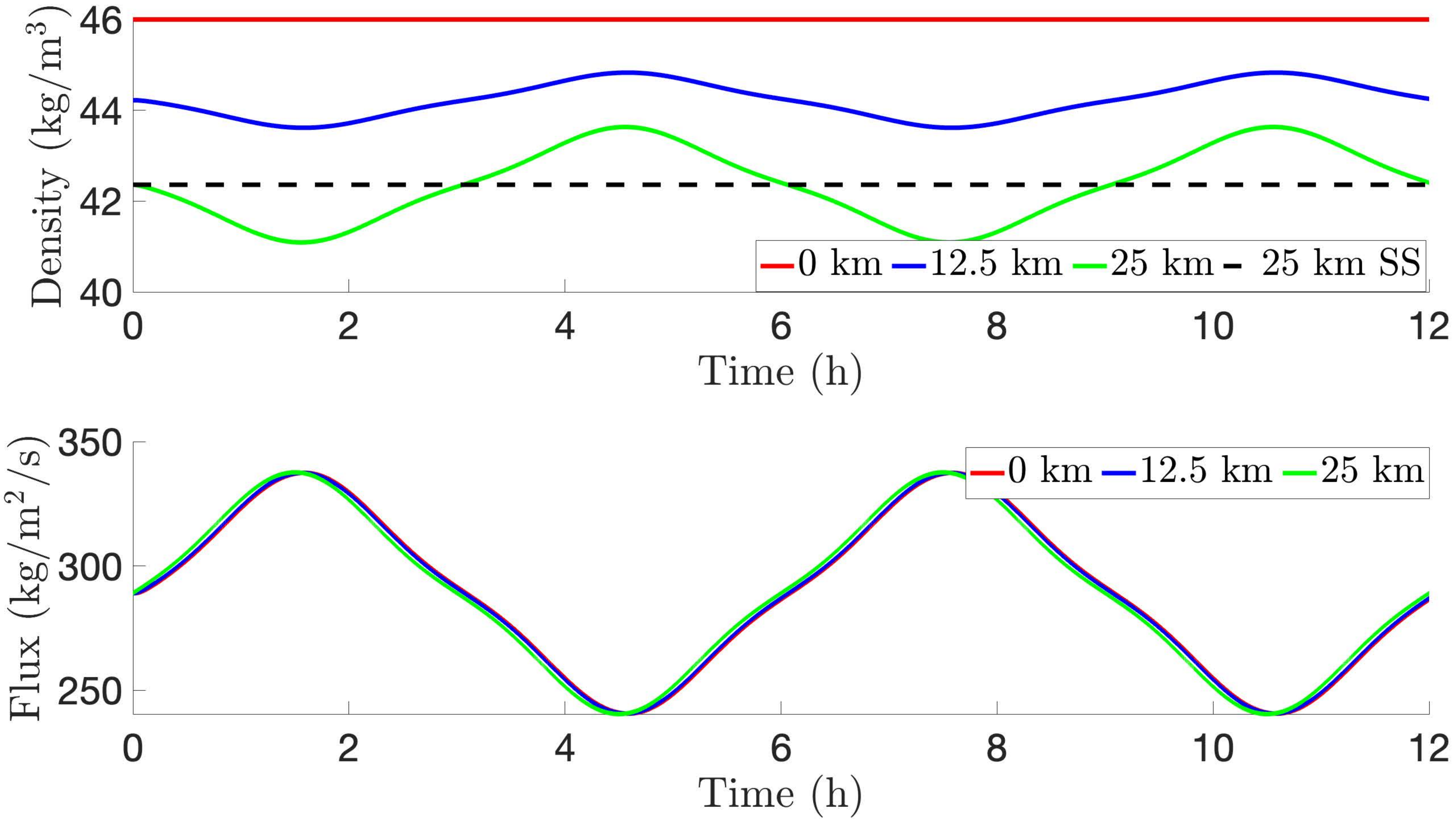}
    \caption{Evolution of density and flow with $\delta U(t)=0$ in the presence of disturbed outlet flow fluctuation.}
    \label{fg6}
\end{figure}

\begin{figure}[h!]
    \centering
    \includegraphics[width=\linewidth]{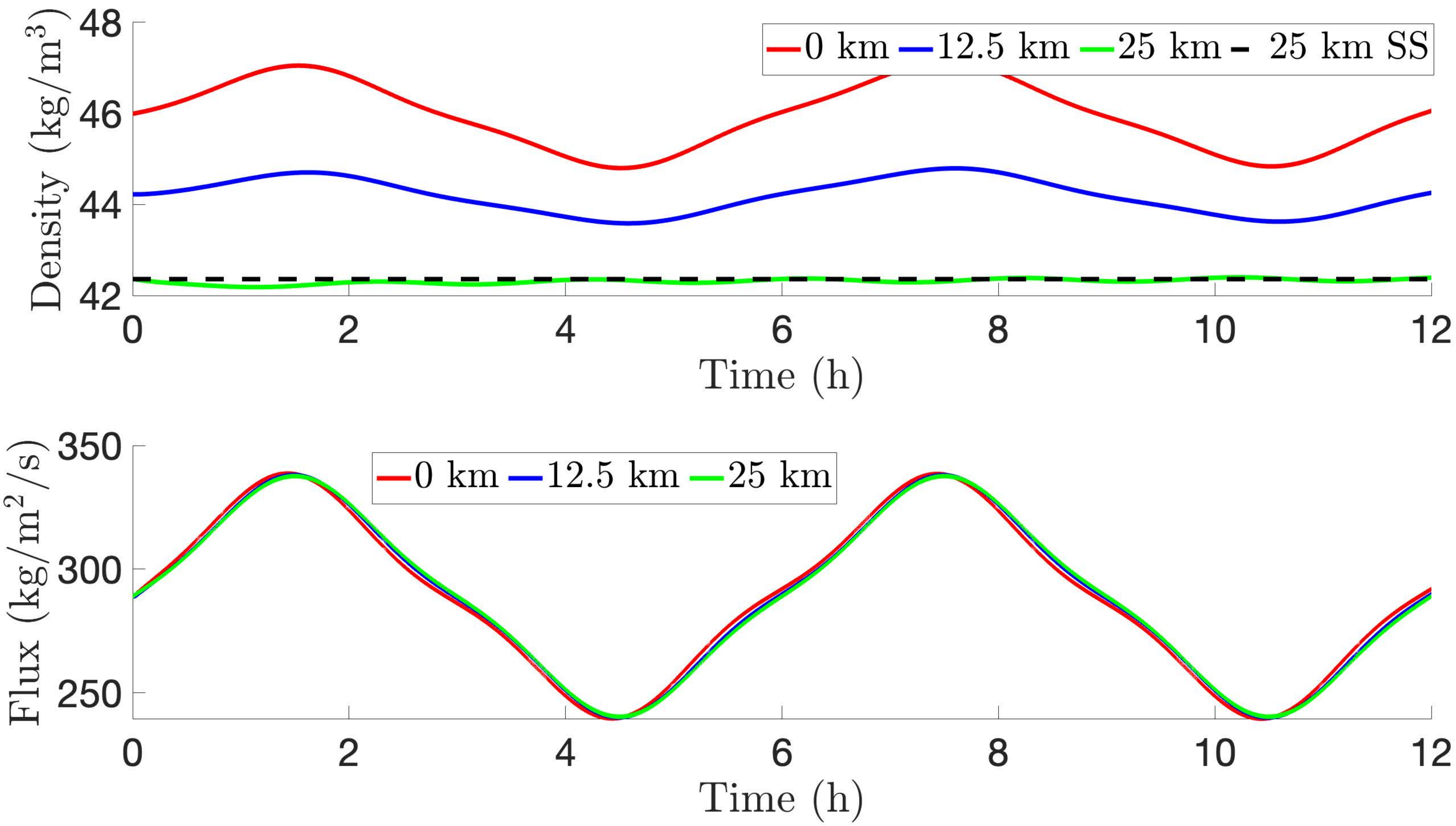}
    \caption{Evolution of density and flow with $\delta U(t)$ given by \eqref{aazmlpsdt1} in the presence of disturbed outlet flow fluctuation.} \vspace{2ex}
    \label{fg7}
\end{figure}

In Figs. \ref{fg6} and \ref{fg7}, we illustrate the evolution of density in the presence of  disturbed outlet flow fluctuation. Fig. \ref{fg6} depicts the scenario where $\delta U(t) = 0$, while Fig. \ref{fg7} shows the case where $\delta U(t)$ is given by \eqref{aazmlpsdt1}. As seen in Fig. \ref{fg7}, the outlet density is regulated to a neighborhood of the steady-state value even in the presence of disturbed outlet flow fluctuations. In contrast, Fig. \ref{fg6} shows that the outlet density fluctuates due to the presence of outlet flow disturbances.


\section{Conclusions}

In this study, we develop an output feedback regulation approach to suppress disturbances from a setpoint in the outlet pressure of gas flowing through a pipeline subject to fluctuating outlet flow. We model the variation in consumption as outlet flow fluctuations generated by a periodic linear dynamic system. The nonlinear isothermal Euler equations including the Darcy-Weisbach friction model are linearized and expressed in canonical form as a coupled $2 \times 2$ hyperbolic PDE. Using an observer-based PDE backstepping controller, we regulate the outlet pressure to a setpoint, even in the presence of outlet flow variations, by manipulating the inlet pressure. Additionally, we extended the controller to handle bounded uncertainties in the outlet flow fluctuations, ensuring regulation of the outlet pressure variation to a neighborhood of the setpoint. The numerical results have validated the effectiveness of the proposed control strategy, and the approach can be extended to gas pipeline networks.

\section{Acknowledgements}

This study was supported by the LDRD program project ``Stochastic Finite Volume Method for Robust Optimization of Nonlinear Flows'' at Los Alamos National Laboratory.  Research is done at Los Alamos National Laboratory under the auspices of the National Nuclear Security Administration of the U.S. Department of Energy under Contract No. 89233218CNA000001. Report No. LA-UR-24-29852.

\bibliographystyle{IEEEtranS}
\bibliography{citation}





\end{document}